\documentclass[reqno, 10pt]{amsart}
\usepackage{color}

\usepackage{graphicx, enumerate}
\usepackage{amssymb, amsmath, amsthm}
\allowdisplaybreaks

\newcommand{\modulo}{\operatorname{mod}}

\newcommand{\corresponds}{\leftrightsquigarrow}

\usepackage{hyperref}
\hypersetup{
    colorlinks=true, 
    linktoc=all,     
    linkcolor=blue,  
    citecolor=blue
}

\newtheorem{Theorem}{Theorem}[section]
\newtheorem{Proposition}[Theorem]{Proposition}
\newtheorem{Lemma}[Theorem]{Lemma}

\newtheorem{Corollary}[Theorem]{Corollary} 
\newtheorem{Definition}[Theorem]{Definition} 
\newtheorem{Conjecture}{Conjecture}

\theoremstyle{definition}
\newtheorem{Example}{Example}[section]
\newtheorem*{Examples}{Examples}

\theoremstyle{remark}

\newtheorem*{Remark}{Remark}
\newtheorem*{Notation}{Notation}

\numberwithin{equation}{section}

\begin{document}

\title{Glaisher's divisors and infinite products}

\author[H.~S.~Bal]{Hartosh Singh Bal}
\address{The Caravan\\
Jhandewalan Extn., New Delhi 110001, India}
\email{hartoshbal@gmail.com}

\author[G.~Bhatnagar]{Gaurav Bhatnagar}
\address{Ashoka University, Sonipat, Haryana, India}
\email{bhatnagarg@gmail.com}
\urladdr{https://www.gbhatnagar.com}

\date{\today}


\maketitle


\begin{abstract}
Ramanujan gave a recurrence relation for the partition function in terms of the sum of the divisor function $\sigma(n)$.  In 1885, J.W. Glaisher considered seven divisor sums closely related to the sum of the divisors function. 
We develop a calculus to associate a generating function with each of these divisor sums. This yields analogues of Ramanujan's recurrence relation for several partition-theoretic functions as well as $r_k(n)$ and $t_k(n)$,  functions counting the number of ways of writing a number as a sum of squares (respectively, triangular) numbers.  As by-products of this association, we obtain several convolutions, recurrences and congruences for divisor functions. We give alternate proofs of two classical theorems, one due to Legendre and the other---Ramanujan's congruence $p(5n+4) \equiv 0 \pmod 5$.
\end{abstract}

\section{Introduction}
Let $\sigma(n)$ be the sum of divisors of $n$ and $p(n)$ the number of unordered integer partitions of $n$. There is a famous recurrence relation connecting them:
\begin{equation}\label{powers1}
np(n) = \sum_{i=1}^{n}\sigma(i)p(n-i).
\end{equation}
This has been found in Ramanujan's work (see \cite[p.~108]{Berndt1994}). 
There are several divisor functions studied in the literature. One can ask whether there are analogues of \eqref{powers1} which relate divisor functions to partition-theoretic functions.

The purpose of this paper is to answer this question for 7 sums over divisors, studied 
by J.\ W.\ Glaisher~\cite{Glaisher1885}. Each of these can 
be expressed in terms of $\sigma(n)$, where we take $\sigma(n)=0$ whenever $n$ is not a positive integer. 
%
%
The seven sums, which we call Glaisher's divisors, are listed in Williams~\cite[p.\ xvi]{Williams2011}, and are as follows.
\begin{subequations}
\begin{align}
 d_1(n) &= \sum_{\substack{d\mid n}}d=\sigma(n); \label{d1-def} \\ 
d_2(n) &=\sum_{\substack{d\mid n\\ d \text{ odd}}} d=\sigma(n) - 2\sigma(n/2); \label{d2-def} \\ 
 d_3(n) &=\sum_{\substack{d\mid n\\d \text{ even}}} d=2\sigma(n/2); \label{d3-def} \\ 
 d_4(n) &=\sum_{\substack{d\mid n\\n/d \text{ odd}}}d=\sigma(n) - \sigma(n/2); \label{d4-def} \\ 
 d_5(n) & =\sum_{\substack{d\mid n\\n/d \text{ even}}}d=\sigma(n/2);\label{d5-def} \\ 
 d_6(n) &=\sum_{\substack{d\mid n}}(-1)^{d-1}d=\sigma(n) - 4\sigma(n/2); \label{d6-def} \\ 
\intertext{and,} 
d_7(n) &= \sum_{\substack{d\mid n}}(-1)^{n/d-1}d=\sigma(n) - 2\sigma(n/2). \label{d7-def}
\end{align}
The expressions of Glaisher's divisors in terms of the divisor function $\sigma(n)$ on the right-hand side of each of the above follow by elementary number-theoretic considerations. 

In this paper, we find connections of Glaisher's divisors with the following functions (see 
Andrews and Eriksson~\cite{AE2004} for an introduction to partition functions).
\begin{align*}
p_o(n) &:=p(n\; | \text{ odd parts}),  \text{ the number of partitions with all parts odd};\\
p_d(n) &:=p(n\; | \text{ distinct parts}), \text{ the number of partitions with distinct parts};\\
p_e(n) &:=p(n\; | \text{ even parts}), \text{ the number of partitions with all parts even};\\
\overline{p}(n) &:=  \text{the number of overpartitions of $n$};
\end{align*}
Of these, only the last one---introduced by Corteel and Lovejoy~\cite{CL2004}---is not self-explanatory. An overpartition is a partition where the first occurrence of a part may be overlined. For example, $5+3+\overline{2}+2+\overline{1}+1+1$ is an overpartition of $15$. 

In addition, two other number-theoretic functions appear in our results.
\begin{align*} 
r_m(n) &:= \text{the number of ways of writing $n$ as an ordered sum of $m$ squares of integers;} \\
t_m(n) &:= \text{the number of ways of writing $n$ as an ordered sum of $m$ triangular numbers.}
\end{align*}

The analogues of \eqref{powers1} we find are all recurrence relations connecting Glaisher's divisors to the above functions. Two other divisor functions appear naturally, as special cases. Let 
$\sigma_3(n)$ be the sum of cubes of the divisors of $n$. Then we require:
\begin{align}
\overline{\sigma}_3(n) := \sum_{\substack{d\mid n\\n/d \text{ odd}}} d^3 = \sigma_3(n) - \sigma_3(n/2);\label{sigmabar-def}\\
\intertext{and} 
\widetilde{\sigma}_3(n) := \sum_{d\mid n}(-1)^{d-1}d^3 = \sigma_3(n) -16\sigma_3(n/2).
\label{sigmacurl-def}
\end{align}
\end{subequations}

Ramanujan had the following very general entry which indicates an approach to recurrences such as \eqref{powers1}. The following is from Berndt~\cite[Entry 12a, p.~28]{Berndt1989}.
Let 
\begin{subequations}
\begin{equation*}
f(q)=\sum_{k=1}^\infty \frac{A_k q^k}{k},
\end{equation*}
and $P_k$ be defined for $k=0, 1, 2, \dots$ by
\begin{equation*}\label{entry2.10.12a-1}
e^{f(q)} =\sum_{n=0}^\infty P_n q^n.
\end{equation*}
Then $P_0=1$ and for $n\ge 1$,
\begin{equation*}\label{entry2.10.12a}
nP_n = \sum_{k=1}^n A_k P_{n-k}.
\end{equation*}
\end{subequations}
Our approach is a mild modification of Ramanujan's approach. We begin with a definition.
 \begin{Definition}[Series-divisor]\label{def-series-divisor}
Let $$A(q)= \sum_{k=0}^{\infty}a(k)q^k$$
be a formal power series. The series divisor $(\sigma^A(k))$ for $k=1, 2, \dots$ is defined by the equation
%
\begin{subequations}
\begin{equation}\label{sigma-divisor-def-a}
\sum_{k=1}^{\infty}ka(k)q^k = \sum_{k=1}^{\infty}\sigma^{A}(k)q^k\sum_{k=0}^{\infty}a(k)q^k.
\end{equation}
 We define
 \begin{equation}\label{sigma-divisor-def-b}
 \sigma^A(r)= 0 \text{ unless } r\in \mathbb{Z}  \text{ and } r>0. 
 \end{equation}
 \end{subequations}
\end{Definition}
 
Alternatively, \eqref{sigma-divisor-def-a} specifies the recurrence
\begin{equation}\label{divisors}
na(n) = \sum_{i=1}^{n}\sigma^{A}(i)a(n-i),
\end{equation} 
which can also be used to define the series-divisor $\sigma^A(n)$, for $n=1, 2, \dots$. 

We develop a calculus to associate a divisor function with an infinite product. In turn, this infinite product is the generating function of the relevant partition-theoretic or arithmetic sequence. The approach is elementary. 
 Nevertheless, we succeed in finding several new results which fit in well with the existing literature.

Note that when the formal power series $A(q)$ is the generating function for the number of partitions of $p(n)$, that is,
$$\prod_{k=0}^{\infty}\frac{1}{(1-q^k)} = \sum_{n=0}^\infty p(n)q^n,$$ 
then \eqref{divisors} reduces to \eqref{powers1}. This motivates the term ``sigma-divisor''.

Evidently, the series-divisor is nothing but the coefficients of the logarithmic derivative of $A(q)$.
Finding the recurrence relation satisfied by a sequence whose generating function is 
known using logarithmic derivatives
is standard in generatingfunctionology (see Wilf~\cite[p.~22]{Wilf2006}). Thus, we usually find the sigma-divisor from a generating function. Here we develop the calculus to go in the reverse 
direction---from the divisor sum to the generating function. 

Our approach was implicit in our earlier work \cite{BB2022a}. In that paper, we applied this idea to 
embed Ramanujan's congruences $p(5n+4) \equiv \tau(5n) \equiv 0  \pmod 5$ into an infinite family of such recurrences. Here $\tau$ is Ramanujan's $\tau$ function \cite[p.~5]{Berndt2006}, defined
by
$$\prod_{k=1}^\infty\bigg( \frac{1}{1-q^k}\bigg)^{24} = \sum_{n=0}^\infty \tau(n+1) q^n.$$

In the current work, the focus is on divisor functions and their correspondence with functions from additive number theory. As an immediate consequence of this correspondence, we obtain recurrence relations connecting the two. In \S\ref{sec:recurrences} we develop the calculus and list these recurrence relations. As special cases, we obtain some more recurrence relations for divisor functions. There are further recurrence relations in \S\ref{sec:powers} which are obtained by using  an elementary idea that
was used very effectively by Gould~\cite{Gould1974}. Next, in 
\S\ref{sec:convolutions}, we round off the applications of our calculus by deriving a pair of formulas for overpartitions of $n$, analogous to results for the partition function due to Euler and Glaisher. 

Aside from \eqref{powers}, there seem to be few results connecting partition-theoretic functions with the multiplicative functions of number theory; see Merca~\cite{Merca2018a} and previous work by the authors~\cite{BB2022a}.  Christopher~\cite{Christopher2019} and Merca~\cite{Merca2021a} have results analogous to Euler's recursion for partitions which is obtained by using the pentagonal number theorem. 
Convolutions of Glaisher's divisors appear in 
\cite{HOSW2002, Hahn2007}. Many such results have been obtained by using Liouville's approach in Williams~\cite{Williams2011}. Regarding congruences for the divisor functions, 
congruences such as $3|\sigma(3n+2)$ and $4|\sigma(4n+3)$ follow from the definition of $\sigma(n)$. Bonciocat ~\cite{bonciocat2003} and Gallardo~\cite{gallardo2007} have proved congruences for the convolution of the sum of divisors functions.
Merca \cite{Merca2021b} has examined congruence sums for $d_2(n)$ over the extended 
pentagonal numbers. Fine~\cite{Fine1988} has shown the application of bilateral $q$-
hypergeometric series evaluations to divisor functions; 
Berdnt~\cite{Berndt2006} uses similar techniques.

The results of this paper complement other work in the area. We illustrate this in \S\ref{sec:new-proofs} by giving alternate proofs of two classical results. One of them is due to Legendre:
$$t_4(j)=\sigma(2j+1).$$
We also give an alternate proof of Ramanujan's famous recurrence
$$p(5n+4)\equiv 0 \pmod 5.$$
Our proof gives an example of an equivalence between a congruence result from partition theory with one involving the convolution of the sum of divisors function $\sigma(n)$. We conclude with another example of this kind involving overpartitions.

\section{Recurrences for the sum of divisors function} \label{sec:recurrences}
We begin our development of a calculus of series-divisors which allows us to virtually read-off from the divisor function the corresponding generating function of an appropriate partition function.

%
 The following is the first of three useful lemmas about series-divisors (recall Definition~\ref{def-series-divisor}).
\begin{Lemma}[Addition Lemma]\label{lemma:addition}
Let $\sigma^{A}(k)$ and $\sigma^{B}(k)$ be the series-divisors for the power series $A(q)$ and 
$B(q)$. Then $\sigma^{A}(k) + \sigma^{B}(k)$ is the series-divisor for $A(q)B(q)$.
That is, the series-divisor of the product of two power series is the sum of the respective series-divisors.
\end{Lemma}
\begin{proof}
The proof is immediate from the fact that the log of a product is the sum of the logs. 
%
\end{proof}
The following proposition allows us to compute series-divisors in many cases of interest.
\begin{Proposition}[Calculus of series-divisors I]\label{calc1}
Let $A(q)$ and $(\sigma^A(k))$ be as in Lemma~\ref{lemma:addition}. For $k=1, 2, \dots, $ we have the following. 
\begin{enumerate}[(i)]
\item The series divisors of $A(q)^{-1}$ are given by $\big(-\sigma^A(k)\big)$. 
\item The series divisors of $A(-q)$ are given by $\big ((-1)^k\sigma^A(k)\big)$. 
\item The series divisors of $A(q^n)$ are given by $\big(n\sigma^A(k/n)\big)$.
\item The series divisors of $A(q)^n$ are given by $\big(n\sigma^A(k)\big)$.
\end{enumerate}
\end{Proposition}
\begin{proof}
While parts (i) and (iv) follows from the addition lemma, the other two parts follow from \eqref{divisors}. For part (iii), recall that we take $\sigma^A(r)=0$ unless $r$ is a positive integer.  
\end{proof}

\begin{Corollary}[Recursion for powers]
Let $a_m(k)$ be the coefficients of the $m$th powers of the power series $A(q)$, that is,
$$A(q)^m= \sum_{k=0}^{\infty}a_m(k)q^k.$$
 Then  
\begin{equation}\label{powers}
na_m(n) = m\sum_{j=1}^{n}\sigma^{A}(i)a_m(n-j).
\end{equation}
\end{Corollary}
\begin{proof}
This is immediate from Proposition~\ref{calc1} (part (iv)) and \eqref{divisors}.
\end{proof}

%

\begin{Examples}[The calculus of series-divisors II] The following will allow us to compute the series-divisor for a host of infinite products. The first is trivial. The rest follow from the previous ones using Proposition~\ref{calc1}.
\begin{enumerate}
\item The series-divisor for $A(q)=1-q$ is $(-1,-1,-1, \dots)$ since
$$-q=(1-q) \sum_{k=1}^\infty (-1) q^k.$$
\item The series-divisor for $(1-q)^{-1}$ is $(1,1,1, \dots).$ 
\item The series-divisor for $1+q$ is $(1,-1,1, \dots).$ 
\item The series-divisor for $(1+q)^{-1}$ is $(-1,1,-1, \dots).$
\item If $A(q)=1-q^n$, then the series-divisor is the sequence 
$$
\sigma^A(k) =
\begin{cases} 
-n &  \text{if } k =mn, \text{ for some } m\in \mathbb{N}\\
0 & \text{otherwise}.
\end{cases}
$$
\item If $A(q)=(1-q^n)^{-1}$, then
$$
\sigma^A(k) =
\begin{cases} 
n &  \text{if } k =mn, \text{ for some } m\in \mathbb{N},\\
0 & \text{otherwise}.
\end{cases}
$$
\item If $A(q)=1+q^n$, then
$$
\sigma^A(k) =
\begin{cases} 
(-1)^{m-1} n &  \text{if } k =mn, \text{ for some } m\in \mathbb{N},\\
0 & \text{otherwise}.
\end{cases}
$$

\item Finally, the series-divisor for $A(q)=(1+q^n)^{-1}$ is
$$
\sigma^A(k) =
\begin{cases} 
(-1)^{m} n &  \text{if } k =mn, \text{ for some } m\in \mathbb{N},\\
0 & \text{otherwise}.
\end{cases}
$$
\end{enumerate}
\end{Examples}

\begin{Remark}[Alternative recursion for the binomial coefficients]
Combining part (1) the above with the recursion for powers, we see that 
the series-divisor for $A(q)=(1-q)^m$ satisfies the recurrence equation 
$$na_m(n) = -m\sum_{j=1}^{n}a_m(n-j) \text{ with } a(0)=1.$$
The first few terms of $a_m(j)$ for $j=0, 1, 2, 3$ are seen to be $1$,$-j$, ${j(j-1)}/{2}$ and $-{j(j-1)(j-2)}/{6}$. From here it it not difficult to obtain the usual formula for the binomial coefficients and prove it satisfies the above recursion. 
\end{Remark}

\begin{Example}[The recurrence \eqref{powers1}]\label{ram-ex1}
Let $p(n)$ be the number of integer partitions of $n$. We apply the addition lemma term by term to the generating function
\begin{equation}\label{partition-gf}
P(q)=\prod_{k=1}^{\infty}\frac{1}{1-q^k} = \sum_{k=0}^\infty p(k)q^k.
\end{equation}
From the above, the series divisors for $(1-q)^{-1}$, $(1-q^2)^{-1},\dots$ are easy to find, and then they are summed using the Addition Lemma. Here we have used  $\leftrightsquigarrow$ to denote the correspondence.
\begin{align*}
\frac{1}{1-q} &\leftrightsquigarrow (1, 1, 1, 1, 1, 1, 1, 1, \dots )\\
\frac{1}{1-q^2} &\leftrightsquigarrow (0, 2, 0, 2, 0, 2, 0, 2, \dots )\\
\frac{1}{1-q^3} &\leftrightsquigarrow (0, 0, 3, 0, 0, 3, 0, 0, \dots )\\
\frac{1}{1-q^4} &\leftrightsquigarrow (0, 0, 0, 4, 0, 0, 0, 4, \dots )\\
\vdots &  \\
\prod_{k=1}^{\infty}\frac{1}{1-q^k} &\leftrightsquigarrow (1, 3, 4, 7, \dots )
\end{align*}
Summing each component, we see that the sigma-divisor is $\sigma(i)$, the sum of divisors function. 
The recurrence \eqref{powers1} is thus a special case of \eqref{divisors}.
\end{Example}

\begin{Notation}
We use the notation $\leftrightsquigarrow$ to indicate the correspondence between sigma-divisors 
$\sigma^A(k)$ and infinite products $A(q)$. 
\end{Notation}

%

The calculations in Example~\ref{ram-ex1} are easily reversed to associate sums of divisors with an infinite product. In the rest of this section we apply this to each one of Glaisher's divisors. The associated infinite products are interesting in their own right, and lead to analogues of \eqref{powers1} with other partition-theoretic functions. 

In the following, the generating functions of the various partition-theoretic functions are required. 
These can be found in \cite{AE2004} and \cite{CL2004}. In addition, we require two 
theta functions from Berndt~\cite[p.~7]{Berndt2006}:
\begin{align}
\varphi(q) &:=\sum_{k=-\infty}^\infty q^{k^2} 
=   \prod_{k=1}^{\infty}{(1 + q^{2k-1})^2}{(1 - q^{2k})} \label{phi}\\
\intertext{ and }
 \psi(q) &:=\sum_{k=0}^\infty q^{\frac{k(k+1)}{2}} = 
 \prod_{k=1}^{\infty}\frac{(1-q^{2k})}{(1-q^{2k-1})} 
 \label{psi}.
 \end{align}
The factorizations into infinite products are due to Gauss and special cases of Jacobi's triple product identity; we refer to Andrews, Askey and Roy~\cite[p.~500]{AAR1999} (a small typo in \cite[(eq. (10.4.8)]{AAR1999} is corrected in \eqref{psi}). Note that these are the only identities
in \cite{AAR1999} where the sums have positive coefficients.

\begin{Theorem}\label{th:assoc} Let $d_1(n), d_2(n), \dots, d_7(n)$ be Glaisher's divisor functions. They are sigma-divisors of the following products, which in turn, are generating functions as given below. 
\begin{subequations}
\begin{align}
d_1(n)=\sigma(n) &\leftrightsquigarrow \prod_{k=1}^{\infty}\frac{1}{1-q^{k}} 
= \sum_{k=0}^\infty p(k)q^k; \label{d1} \\
d_2(n) &\leftrightsquigarrow \prod_{k=1}^{\infty}\frac{1}{1-q^{2k-1}}
= \sum_{k=0}^\infty p_o(k)q^k ; \label{d2} \\
d_3(n) = 2d_5(n) &\leftrightsquigarrow \prod_{k=1}^{\infty}\frac{1}{1-q^{2k}}
= \sum_{k=0}^\infty p_e(k)q^k ; \label{d3}\\
2d_4(n) &\leftrightsquigarrow  \prod_{k=1}^{\infty}\frac{1+q^k}{1-q^k}
= \sum_{k=0}^\infty \overline{p}(k)q^k ; \label{d4-a}\\
(-1)^{n+1}2d_4(n) &\leftrightsquigarrow \prod_{k=1}^{\infty}(1+q^{2k-1})^2(1-q^{2k}) 
=\varphi(q) =\sum_{k=-\infty}^\infty q^{k^2} ; \label{d4-b}\\
d_6(n) &\leftrightsquigarrow    \prod_{k=1}^{\infty}\frac{(1-q^{2k})}{(1-q^{2k-1})} = \psi(q)
=\sum_{k=0}^\infty q^{\frac{k(k+1)}{2}}; \label{d6} \\
d_7(n) &\leftrightsquigarrow \prod_{k=1}^{\infty}(1+q^k)
= \sum_{k=0}^\infty p_d(k)q^k .  \label{d7}
\end{align}
\end{subequations}
\end{Theorem}
\begin{proof}
The proof of \eqref{d1} outlined in Example~\ref{ram-ex1} extends easily to prove \eqref{d2} and \eqref{d3}. Here only the odd (respectively, even) divisors are there, with only the corresponding terms in the infinite product. That $d_3(n)=2d_5(n)$ follows from the right-hand side of \eqref{d3-def}
and \eqref{d5-def}.

Next we consider the sequence $\big( d_7(n)\big)$. The sequence is generated by summing $(1, -1, 1, -1, \dots,)$, 
$(0, 2, 0, -2, \dots,)$ and so on. Thus
$$d_7(n) \leftrightsquigarrow \prod_{k=1}^{\infty}(1+q^k)$$
which is the generating function of partitions with distinct parts. This gives \eqref{d7}. 

From here we obtain 
\begin{align*} 2d_4(n) &= d_1(n)+ d_2(n) 
\corresponds \prod_{k=1}^{\infty}\frac{1}{1-q^{k}}  \prod_{k=1}^{\infty}\frac{1}{1-q^{2k-1}} 
= \prod_{k=1}^{\infty}\frac{1+q^k}{1-q^k},
\end{align*}
where we have used the analytic form of Euler's theorem that
the number of partitions with distinct parts equals the number of partitions with odd parts. This shows 
\eqref{d4-a}.

To show \eqref{d4-b}, consider the following:
\begin{align*}
\varphi(q)= 
\prod_{k=1}^{\infty} & (1-(-q)^{2k-1})^2(1-(-q)^{2k}) 
\corresponds -(-1)^n 2 d_2(n) -(-1)^n d_3(n)  \\
&= (-1)^{n+1}\Big(\big( 2(\sigma(n)-2\sigma(n/2)\big)+2\sigma(n/2)\big) = 
(-1)^{n+1} 2 d_4(n).
\end{align*}

Finally, observe that 
\begin{align*}
d_6(n) = \sigma(n)-4\sigma(n/2) = 
d_2(n)-d_3(n)
\corresponds 
\prod_{k=1}^{\infty}\frac{(1-q^{2k})}{(1-q^{2k-1})} = \psi(q).
\end{align*}
\end{proof}

\begin{Remark}
Here is an  interesting application of Euler's ODD $=$ DISTINCT theorem. Note that
\begin{align*}
d_7(n) \leftrightsquigarrow & \prod_{k=1}^{\infty}(1+q^k) 
 = \prod_{k=1}^{\infty}\frac{1}{1-q^{2k-1}} 
 \leftrightsquigarrow d_2(n).
 \end{align*}
Thus $$d_7(n)=d_2(n)$$ and 
both equal $$d_1(n)-d_3(n)=\sigma(n)-2\sigma(n/2).$$
This is an alternate partition-theoretic proof of the right hand side of \eqref{d7-def}. 
\end{Remark}

From the correspondence between Glaisher's divisor functions and partition-theoretic generating functions, we obtain analogues of Ramanujan's recurrence \eqref{powers1}.
\begin{Theorem} The following recurrence relations hold. 
\begin{subequations}
\begin{align}
np_e(n) &=\sum_{i=1}^{n}d_3(i)p_e(n-i)= 2\sum_{i=1}^{n}\sigma(i/2\big)p_e(n-i) ;\\
np_d(n) &=\sum_{i=1}^{n}d_2(i)p_{o}(n-i)= \sum_{i=1}^{n}d_7(i)p_d(n-i) ; \label{pd-rec-a} \\
&=\sum_{i=1}^{n}\big(\sigma(i) - 2\sigma(i/2)\big)p_d(n-i); \label{pd-rec-b} \\
n\overline{p}(n) &= 2\sum_{i=1}^{n}d_4(i)\overline{p}(n-i)= 2\sum_{i=1}^{n}\big(\sigma(i) - \sigma(i/2)\big)\overline{p}(n-i).
\end{align}
\end{subequations}
\end{Theorem}
\begin{proof}
These recurrence relations are just \eqref{divisors} when applied to 
\eqref{d3}, \eqref{d2} and \eqref{d7}, and  \eqref{d4-a}, respectively. In \eqref{pd-rec-a} and 
\eqref{pd-rec-b}, we use Euler's theorem $p_d(n)=p_o(n)$.
\end{proof}

\begin{Theorem} For integers $n\geq0$ and $m>0$, let $r_m(n)$, represent the number of ways $n$ can be written as an ordered sum of $m$ squares, and let $t_m(n)$, represent the number of ways $n$ can be written as an ordered sum of $m$ triangular numbers. Then
\begin{subequations}
\begin{align}
nr_m(n) &= m\sum_{j=1}^{n}2 (-1)^{j+1}\big(\sigma(j) - \sigma(j/2)\big)r_m(n-j); \label{squares} \\
\intertext{and,}
nt_m(n) &= m\sum_{j=1}^{n}\big(\sigma(j) - 4\sigma(j/2)\big)t_m(n-j) \label{triangular}.
\end{align} 
\end{subequations}
\end{Theorem}

\begin{proof}
Since the coefficient of $q^n$ of $\varphi(q)^m$ is the number of ways $n$ can be written as a sum of $m$ squares, by \eqref{powers} we have \eqref{squares}. Similarly, the coefficients of
$\psi(q)^m$ are the number of ways $n$ can be written as a sum of $m$ triangular numbers;  by \eqref{powers} we have \eqref{triangular}.
\end{proof}

\begin{Remark} The calculus of series-divisor presented above is implicit in \cite{BB2022a}, where we consider any power series generated by an infinite product of the form 
$$\prod_{k=0}^{\infty}(1-q^k)^{f(k)}.$$ In \cite{BB2022a}, $f(k)$ can take complex values or could even be a polynomial with complex coefficients. Here we only require integral powers, so we provided a simpler exposition.
\end{Remark}  

Next we take special cases of \eqref{squares} and \eqref{triangular} to obtain several convolution recurrences for Glaisher's divisor functions in terms of $\sigma(n)$, $\overline{\sigma}_3(n)$ and
$\widetilde{\sigma}_3(n)$ (see \eqref{sigmabar-def} and \eqref{sigmacurl-def}). 

To take special cases, we use the following well-known identities given in \cite[Theorems 3.3.1 and 3.3.4]{Berndt2006}.
For $n>0$, we have
\begin{subequations}
\begin{align}\label{ntsquare-a}
r_1(n) &=
\begin{cases}
2, & \text{if } n=j^2 \text{ for some } j>0; \\
0, &\text{otherwise}
\end{cases}\\
\label{ntsquare-b}
r_4(n) &= 8\big(\sigma(n)- 4\sigma(n/4)\big) \\
\label{ntsquare-c}
r_8(n) &=(-1)^{n+1}16\big(\sigma_3(n)-16\sigma_3(n/2)\big)=(-1)^{n+1}16\widetilde{\sigma}_3(n).
\end{align}
\end{subequations}
where $r_1(0)=r_4(0)=r_8(0)=1$. 

Next, we use the following for the sum of triangular numbers (see~\cite[Theorems 3 and 5]{ORW1995}):
\begin{subequations}
\begin{align}\label{nttriangle-a}
t_1(n) &=
\begin{cases}
1, & \text{if } n=j(j+1)/2 \text{ for some } j>0; \\
0, &\text{otherwise}
\end{cases}\\
\label{nttriangle-b}
t_4(n) &= \sigma(2n+1)\\
\label{nttriangle-c}
t_8(n) &= \sigma_3(n+1) - \sigma_3((n+1)/2)= \overline{\sigma}_3(n+1)
\end{align}
\end{subequations}
where $t_1(0)=t_4(0)=t_8(0)=1$. 

\begin{Theorem}\label{theorem1}
Let $n$ be a positive integer. Then we have the following recursions:
\begin{align}
\label{div2}
d_4(n)+2\sum_{j=1}^\infty (-1)^{j}d_4(n-j^2) &= 
\begin{cases}
(-1)^{n-1}n, & \text{if } n=j^2 \text{ for some } j; \\
0, &\text{otherwise.}
\end{cases} 
\\
\label{div1}
\sum_{j=0}^\infty d_6\big(n-\frac{j(j+1)}{2}\big) &= 
\begin{cases}
n, & \text{if } n=j(j+1)/2 \text{ for some } j; \\
0, &\text{otherwise.}
\end{cases}
 \\
\label{div4}
8\sum_{j=1}^{n-1} (-1)^{j+1} d_4(j)\big(\sigma(n-j)- 4\sigma(\frac{n-j}{4})\big)
& =n\big(\sigma(n)- 4\sigma(\frac{n}{4})\big)+(-1)^{n}d_4(n) \\
\label{div3}
4\sum_{j=1}^\infty d_6(j)\sigma(2n+1-2j) &=n\sigma(2n+1).
\\
 \label{div6}
16\sum_{j=1}^{n-1}d_4(j)\widetilde{\sigma}_3(n-j) &= d_4(n) - n\widetilde{\sigma}_3(n).\\
\label{div5}
8\sum_{j=1}^{n} d_6(j)\overline{\sigma}_3(n+1-j) &= n\overline{\sigma}_3(n+1) .
\end{align}
\end{Theorem}
\begin{Remark}
Note that \eqref{div2} is due to Liouville; see \cite[Theorem 6.1]{Williams2011}. Williams~\cite{Williams2011} also uses the notation $\sigma^*(n)$ for $d_4(n)$ and 
$\widetilde{\sigma}(n)$ for $d_6(n)$. 
\end{Remark}
\begin{proof}
These results are special cases of \eqref{squares} and \eqref{triangular}.
Plug in $m=1,4$ and $8$ in \eqref{squares}, and use \eqref{ntsquare-a}-\eqref{ntsquare-c}, to obtain \eqref{div2}, \eqref{div4} and \eqref{div6} (respectively). 
Similarly, \eqref{triangular} and \eqref{nttriangle-a}-\eqref{nttriangle-c} yield \eqref{div1}, \eqref{div3} and \eqref{div5}.
\end{proof}
\begin{Remark} We can obtain analogous results by specializing \eqref{squares} and \eqref{triangular} using other values of  $m$. For example, see \cite[Theorems 
3.2.1 and 3.4.1]{Berndt2006} for expressions for $r_m(n)$ with $m=2$ and $6$,  
and  \cite[\S 3]{ORW1995} for expressions for $t_m(n)$ with several small values fo $m$ not considered here. All these give rise to analogous results. 
 This remark applies to the theorems in \S\ref{sec:powers} too. 
\end{Remark}

\section{Powers of a power series} \label{sec:powers}
The next set of results require an old trick involving powers of generating functions.
It is a recurrence relation connecting powers of a generating function. 
%
%
\begin{Lemma}[Power Recursion Lemma]
For any power series $A(q)$, the coefficients $a_u$, $a_v$ of its powers $\big(A(q)\big)^u$ and $\big(A(q)\big)^v$,  where $u$ and $v$ are any two non-zero integers, satisfy: 
\begin{equation}\label{general}
\sum_{k=0}^{n}\big(n-(u/v + 1) k\big)a_u(n-k)a_v(k)=0.
\end{equation}
\end{Lemma}
\begin{Remark} This result hold as long as $A(q)^u$ and $A(q)^v$ are formal power series; see \cite[Prop,~6.1]{BB2022a}. Gould \cite{Gould1974} credits this to Rothe (1793). For the sake of completeness, we repeat the proof from \cite{BB2022a} for integral $u$ and $v$.
\end{Remark}
\begin{proof}
We have $$\sum_{k=0}^{\infty}ka_u(k)q^k=u\big(\sum_{k=1}^{\infty}\sigma^{A}(k)q^k\big)(A(q))^u$$ and 
$$\sum_{k=0}^{\infty}ka_v(k)q^k=v\big(\sum_{k=1}^{\infty}\sigma^{A}(k)q^k\big)(A(q))^u.$$
From here, we obtain
$$\big(\sum_{k=0}^{\infty}ka_u(k)q^k\big)\big(\sum_{k=0}^{\infty}a_v(k)q^k\big) =\frac{u}{v}(\sum_{k=0}^{\infty}ka_v(k)q^k)(\sum_{k=1}^{\infty}a_u(k)q^k).$$
Carrying out the Cauchy product and collecting terms we obtain \eqref{general}.
\end{proof}

Next, as in \cite{BB2022a}, we take $u=1$, $v=m$, with $A(q)=\varphi(q), \psi(q)$ to obtain
\begin{align}\label{square2}
nr_m(n) &=-2\sum_{k=1}^{\infty}\big(n-(m+1)k^2\big)r_m(n-k^2) 
\intertext{and}
\label{triangular2}
nt_m(n)&=-\sum_{k=1}^{\infty}\big(n-(m+1)k(k+1)/2\big)t_m\big(n-k(k+1)/2\big).
\end{align}
Of these, \eqref{square2} appears in Venkov~\cite[p.~204]{Venkov1970} and Williams~\cite[p.~44]{Williams2011}. 
%
%
%
These results yield recurrence relations for $\sigma(n)$, $\overline{\sigma}_3(n)$ and $\widetilde{\sigma}_3(n)$.

\begin{Theorem}\label{theorem3}
Let $n$ be a positive integer. Then we have the following recursions:
\begin{equation}\label{div7}
2n\sigma(2n+1)=\sum_{j=1}^\infty \big(5j(j+1)- 2n \big) \sigma\big(2n+1-j(j+1)\big),
\end{equation}
\begin{equation}\label{div8}
n\big(\sigma(n) - 4\sigma(\frac{n}{4})\big)- 2\sum_{j=1}^\infty (5j^2 - n)\big(\sigma(n -j^2)-4\sigma(\frac{n -j^2}{4})\big)=
\begin{cases}
n, & \text{if } n=j^2 \text{ for some } j; \\
0, &\text{otherwise}
\end{cases}
\end{equation}
\begin{equation}\label{div9}
n\overline{\sigma}_3(n+1)=\sum_{j=1}^\infty\big(9\frac{j(j+1)}{2} - n\big)\overline{\sigma}_3(n+1 - \frac{j(j+1)}{2}),
\end{equation}
\begin{equation}\label{div10}
n\widetilde{\sigma}_3(n)-2\sum_{j=1}^\infty (-1)^{j}(9j^2 - n)\widetilde{\sigma}_3(n -j^2)=
\begin{cases}
(-1)^{n-1}n, & \text{if } n=j^2 \text{ for some } j; \\
0, &\text{otherwise}
\end{cases}
\end{equation}
\begin{equation}\label{div11}
\sum_{j=0}^\infty(3j-n)\sigma(2j+1)\overline{\sigma}_3(n+1-j)=0,
\end{equation}
and
\begin{equation}\label{div12}
n(\widetilde{\sigma}_3(n)+(-1)^{n}\big(\sigma(n)-4\sigma(\frac{n}{4})\big)=8\sum_{j=1}^{n-1}(3j- n)(-1)^j\big(\sigma(j)-4\sigma(\frac{j}{4})\big)\widetilde{\sigma}_3(n-j).
\end{equation}
\end{Theorem}
\begin{proof}
These results follow from \eqref{square2} and \eqref{triangular2} by taking $n=1, 4, 8$, and using
\eqref{ntsquare-a}-\eqref{ntsquare-c} and \eqref{nttriangle-a}-\eqref{nttriangle-c}.
\end{proof}

The recurrences in Theorem~\ref{theorem1} and Theorem~\ref{theorem3} yield several congruences for divisor functions.
\begin{Theorem}\label{theorem4}
Let $n$ be a positive integer. Then we have the following congruences:
If $5 \nmid n$, then
\begin{equation}\label{div13}
\sum_{j=0}^\infty \sigma\big(2n+1-j(j+1)\big)\equiv0 \pmod 5.
\end{equation}
%
\begin{equation}\label{div14}
\sum_{j=0}^\infty \sigma\big(n-j(j+1)/2\big) \equiv 
\begin{cases}
0 \; (\modulo 4) & \text{ if } n \neq \frac{k(k+1)}2 \text{ for some }k\\
n  \; (\modulo 4) & \text{ if } n = \frac{k(k+1)}2 \text{ for some } k.
\end{cases}
\end{equation}
If $5 \nmid n$ and $n=4m+3$ or $4m+2$, then
\begin{equation}\label{div15}
\sigma(n) + 2\sum_{j=1}^\infty \sigma(n-j^2)\equiv0 \pmod 5.
\end{equation}
If $(n,3) = 1$, then
\begin{equation}\label{div16}
\sum_{j=0}^\infty \overline{\sigma}_3\Big(n+1 - \frac{j(j+1)}{2}\Big) \equiv0 \pmod {3^2}.
\end{equation}
\end{Theorem}

\begin{proof}
If $5 \nmid n$, consider the identity \eqref{div7} mod $5$ and cancel the $n$ to 
obtain \eqref{div13}.

The identity \eqref{div1} mod $4$ immediately yields
 \eqref{div14}. 
 
 Next, consider identity \eqref{div8} mod $5$ and note that $j^2$ is always of the form $4k$ or $4k+1$ for some integer $k$ so the terms $\sigma((4m+3 - j^2)/4) = 0$. 
Considering the resulting identity mod $5$ gives \eqref{div15}.

Finally,  \eqref{div16} follows from \eqref{div9}.
\end{proof}

Regarding \eqref{div13}, more may be true. Computer experiments suggest the following conjecture.
\begin{Conjecture} 
\begin{equation}\label{conj1}
\sum_{j=0}^\infty \sigma\big(2n+1-j(j+1)\big) \equiv
\begin{cases}
0 \; (\modulo 5) & \text{ if } n \neq \frac{5k(k+1)}2 \text{ for some }k\\
1  \; (\modulo 5) & \text{ if } n = \frac{5k(k+1)}2 \text{ for some } k 
\end{cases}
\end{equation}
\end{Conjecture}

\begin{Remark} Here is another nice trick. Consider the product
$$\psi(q)^5 = \psi(q)^4\psi(q),$$
and recall that $t_4(n)= \sigma(2n+1)$. On comparing coefficients, we obtain
%
$$t_5(n) 
=\sum_{j=0}^\infty \sigma\big(2n+1-j(j+1)\big).$$
By \eqref{triangular} this is a multiple of $5$ when $5 \nmid n$. This gives an alternate proof of 
\eqref{div14}.

A similar approach works for \eqref{div15} and \eqref{div16}.

This process can be iterated and we note that 
$$t_6(n)= \sum_{j,k=0}^{\infty} \sigma(2n+1-j(j+1) - k(k+1))$$ 
is divisible by $6$ if $n$ is relatively prime to $6$ and so on.
\end{Remark}

\section{Convolutions of series divisors} \label{sec:convolutions}
In this section, we give a pair of formulas for $\overline{p}(n)$, the number of overpartitions of $n$. 
%
%
These are analogous to the following two results, due to Euler
~\cite{LE1760-243}
and Glaisher~\cite{Glaisher1891} (respectively).
\begin{equation}\label{euler}
\sigma(n)=\sum_{i=-\infty}^\infty (-1)^{i}\big(n-\frac{i(3i-1)}{2}\big)p\big(n-\frac{i(3i-1)}{2}\big)
\end{equation}
and 
\begin{equation}\label{euler*}
\sum_{i=1}^{n}\sigma(n-i)\sigma(i)=\sum_{i=-\infty}^\infty (-1)^{i+1}\Big(\frac{i(3i-1)}{2}\Big)\Big(n-\frac{i(3i-1)}{2}\Big)p\Big(n-\frac{i(3i-1)}{2}\Big).
\end{equation}
We had derived \eqref{euler} and \eqref{euler*} in \cite{BB2022a}, but had inadvertently omitted the reference to Glaisher. The following lemma formalizes our calculation in \cite{BB2022a}.



%
\begin{Lemma}[Convolution lemma]
Let $A(q)$ be a formal power series with coefficients $a_n$, and let $B(q)=1/A(q)$ have coefficients $b_n$. Then
\begin{equation}\label{sigma}
\sigma^{A}(n)= \sum_{i=0}^{n}(n-i)a(n-i)b(i);
\end{equation}    
and
\begin{equation}\label{sigma*}
\sum_{i=1}^{n}\sigma^{A}(n-i)\sigma^{A}(i) = -\sum_{i=1}^{n}i(n-i)a(n-i)b(i).
\end{equation} 
\end{Lemma}
\begin{proof}
Since 
$$A(q)\sum_{i=1}^\infty\sigma^{A}(i)q^i=\sum_{i=1}^\infty ia(i)q^i,$$
we get $$\sum_{i=1}^\infty\sigma^{A}(i)q^i=B(q)\sum_{i=1}^\infty a(i)q^i$$
and carrying out the Cauchy product gives us \eqref{sigma}.

Since the series divisor of $1/A(q)$ is $-\sigma^A$, we also obtain 
$$-\frac{A(q)}{A(q)}\sum_{i=1}^\infty\sigma^{A}(i)q^i\sum_{i=1}^\infty\sigma^{A}(i)q^i=\sum_{i=1}^\infty ia(i)q^i\sum_{i=1}^\infty ib(i)q^i$$ and carrying out the Cauchy product yields \eqref{sigma*}.
\end{proof}
The results cited above follow from \eqref{sigma} and \eqref{sigma*} by taking $A(q)$ to be the generating function for partitions. The expansion of $A(q)^{-1}$ as a power series is Euler's pentagonal number theorem \cite[p.~500]{AAR1999}:
  $$\prod_{k=1}^{\infty}(1-q^k) = \sum_{n=-\infty}^\infty (-1)^{i}q^{\frac{n(3n-1)}{2}}.$$
The key idea is that the power series expansion of the reciprocal should be available. 
\begin{Theorem}
Let $\overline{p}(n)$ be the number of overpartitions of n. Then
\begin{equation}\label{euler1}
2d_4(n)= 2\big(\sigma(n)-\sigma(n/2)\big)= n\overline{p}(n) + 2\sum_{i=1}^\infty (-1)^{i}(n-i^2)\overline{p}(n-i^2);
\end{equation}
and,
\begin{equation}\label{euler1*}
2\sum_{i=1}^{n}d_4(n-i)d_4(i)= \sum_{i=1}^\infty (-1)^{i+1}i^2(n-i^2)\overline{p}(n-i^2).
\end{equation}
\end{Theorem}
\begin{proof}
By Theorem~\ref{th:assoc}, the series-divisor $2d_4(n)$ is associated with the product 
$$\prod_{k=1}^{\infty}\frac{(1+q^k)}{(1-q^k)} = \sum_{n=1}^\infty \overline{p}(n)q^n$$ 
and 
$$\prod_{k=1}^{\infty}\frac{(1-q^k)}{(1+q^k)} = 1 + 2\sum_{n=1}^\infty (-1)^{n}q^{n^2}.$$
The results follows from \eqref{sigma} and \eqref{sigma*}.
The identity  for the reciprocal of the generating function of $\overline{p}(n)$ follows from 
\eqref{phi} by replacing $q$ by $-q$ and manipulating the infinite products. 
\end{proof}

\section{Two classical results}\label{sec:new-proofs}
In this section we give alternate proofs of two classical theorems, due to Legendre and Ramanujan.
These illustrate the application of our approach to prove both results involving the divisor functions and those in the theory of partitions. 

First we use \eqref{triangular} along with two results due to Melfi~\cite{Melfi1998} to prove \eqref{nttriangle-b}, a result of Legendre.
\begin{Proposition}[Legendre] For $j=0, 1, 2, 3, \dots$, 
$$t_4(j)=\sigma(2j+1).$$
\end{Proposition}
\begin{proof}
By \eqref{triangular} $$jt_4(j) = 4\sum_{i=1}^{j}\big(\sigma(i) - 4\sigma(i/2)\big)t_4(j-i).$$
We show that $\sigma(2j+1)$ satisfies the same recurrence by invoking two convolution identities proved by elementary means by Huard, Ou, Spearman and Williams in~\cite{HOSW2002}:  
$$\sum_{i=1}^{n}\sigma(i)\sigma(2n+1-2i)= \frac{1}{24}\big(2\sigma_3(2n+1)+ (1-3(2n+1)\big)\sigma(2n+1),$$ 
$$\sum_{i=1}^{\lfloor{n/2}\rfloor}\sigma(i)\sigma(2n+1-4i)= \frac{1}{48}\Big(\sigma_3(2n+1)+ \big(2-3(2n+1)\big)\Big)\sigma(2n+1).$$
These are special cases of results of Melfi~\cite{Melfi1998}; see \cite[Th.\ 2]{HOSW2002} and \cite[Th.\ 4)]{HOSW2002}.  
Plugging in these identities, yields the desired result:
$$4\sum_{i=1}^{j}\big(\sigma(i) - 4\sigma(i/2)\big)\sigma(2j+1-2i)=j\sigma(2j+1).$$  
\end{proof}

For more proofs of Legendre's result, see \cite{ORW1995}, \cite{HOSW2002},  \cite[p.~72]{Berndt2006}  and 
\cite{Hahn2007}.


\begin{Theorem}\label{th:equiv1} Let $p(n)$ denote the number of partitions of $n$. Then the following statements are equivalent.
\begin{subequations}
\begin{align}
p(5m+4) &\equiv 0 \pmod 5 \label{ram1} \\
\sum_{i\ge 1}\sigma(i)\sigma(5m+1-i) &\equiv 0 \pmod 5. \label{ram2}
\end{align}
\end{subequations}
\end{Theorem}
\begin{proof}
The key idea is to 
consider \eqref{euler*}  mod $5$, with $n=5m+1$. The left hand side is the convolution sum 
in \eqref{ram2}. As for the right-hand side of \eqref{euler*}, note that
for $i=0, 1, 2 \pmod 5$, the product 
$$\Big(\frac{i(3i-1)}{2}\Big) \Big(n- \frac{i(3i-1)}{2}\Big)\equiv 0\pmod 5.$$
So the only terms that survive in the product are when $i=3, 4$ when
$i(3i-1)/2 \equiv 2 \pmod 5$ and the product is of the form
$2 (5m-1) \equiv -2 \pmod 5$. 
Thus we see that when $n=5m+1$, the sum on the right-hand side of \eqref{euler*} is over terms of the form
$$ (*)p(5 m +1 - k ), \text{ with } k \equiv 2 \pmod 5.$$

To be precise, we 
change the index by replacing $i=5j-1$ and $i=5j-2$, and write  \eqref{euler*} $\pmod 5$ as
\begin{multline}\label{euler*mod5}
\sum_{i\ge 1} \sigma(5m+1-i)\sigma(i) =
\sum_{j=-\infty}^\infty (-1)^{j} 2\bigg(
p\Big(5m+1-\frac{(5j-1)(15j-4)}{2}\Big) 
+ \\
p\Big(5m+1-\frac{(5j-2)(15j-7)}{2}\Big) \bigg)
 \pmod 5.
 \end{multline}

Replacing $m$ by $m+1$, we see that all the terms have $p(5k+4)$ for some $k$. If they are all $0$ by \eqref{ram1}, we immediately obtain \eqref{ram2}.

Conversely, we assume \eqref{ram2}, and
use induction to show that $p(5m+4)\equiv 0 \pmod 5$. For $m=0$, the result is true
because $p(4)=5$.  Suppose it is true for numbers less than $m$. 

Note that at $j=0$, $(5j-1)(15j-4)/{2}$ is $2$. 
Thus using the right-hand side of \eqref{euler*mod5} with $m$ replaced by  $m+1$, we can write $p(5m+6-2)=p(5m+4)$ in terms of $p(5k+4)$ with $k<m$. 
Since the sum is $0 \pmod 5$ by 
\eqref{ram2}, we obtain \eqref{ram1} by induction.
\end{proof}

%
Next, we give a new proof of one of Ramanujan's congruences for partitions. This proof uses a result of Jacobi, as well as \eqref{euler*} (which relies on Euler's pentagonal number theorem). Both these are special cases of Jacobi's triple product identity. 

\begin{Proposition}[Ramanujan]
Let $p(n)$ be the number of partitions of $n$. Then
$$p(5m+4)\equiv 0 \pmod 5.$$
\end{Proposition}
\begin{proof} In view of Theorem~\ref{th:equiv1}, it is enough to show \eqref{ram2}.
We use the following result due to Jacobi~\cite[p.~500]{AAR1999}:
\begin{equation}\label{jacobi1}
\prod_{k=1}^\infty (1-q^k)^3 =
 \sum_{k=0}^\infty (-1)^k (2k+1) q^{\frac{k(k+1)}{2}}.
\end{equation}
Let $p_3(n)$ be defined from the generating function
$$A(q)= \prod_{k=1}^\infty \frac{1}{(1-q^k)^3} =\sum_{n=0}^\infty p_3(n)q^n.$$
(The quantity $p_3(n)$ is the number of partitions of $n$ where each part can occur in $3$ colors.)
Then \eqref{sigma*} and Proposition~\ref{calc1} (part (iv)) yield
$$9 
\sum_{i\ge 1}\sigma(i)\sigma(n-i) 
=  \sum_{k=0}^\infty (-1)^{k+1} \frac{k(k+1)(2k+1)}{2}\Big(n - \frac{k(k+1}2\Big) 
p_3\Big(n-\frac{k(k+1)}2\Big) .
$$
Now when $n=5m+1$, consider each term for $k=0,1,2,3,4$. Clearly, each term is
$0\pmod 5$. 
  This implies \eqref{ram2}. 
\end{proof}

We conclude with a theorem for overpartitions on the lines of Theorem~\ref{th:equiv1}.
\begin{Theorem}\label{th:equiv2} Let $\overline{p}(n)$ denote the number of partitions of $n$. Then the following statements are equivalent.
\begin{subequations}
\begin{align}
\overline{p}(4m+3) &\equiv 0 \pmod 8 \label{overp1} \\
\sum_{i\ge 1} d_4(i)d_4(4m-i)  &\equiv 0 \pmod 4. \label{overp2}
\end{align}
\end{subequations}
\end{Theorem}

\begin{proof}
We first show \eqref{overp1} implies \eqref{overp2}. We consider \eqref{euler1*} when $n=4m$. For all $i$, $i^2\equiv 0 \text{ or } 1 \pmod 4$, so 
$$4m-i^2\equiv 0 \text{ or } 3 \pmod 4.$$
When $4m-i^2\equiv 0\pmod 4$, then $i^2$ and $4m-i^2$ are both multiples of $4$, so 
we see that the left hand side of \eqref{euler1*} is divisible by $16$. The 
convolution
$$\sum_{i\ge 1} d_4(i)d_4(4m-i)$$ is divisible by $8$, and so by $4$.  
When $4m-i^2\equiv 3\pmod 4$, then by \eqref{overp1} the right hand side of \eqref{euler1*} is divisible by $8$ and again \eqref{overp2} holds. 

Conversely, if \eqref{overp2} holds, only the terms when $4m-i^2\equiv 3 \pmod 4$ survive in the sum
$$ \sum_{i=1}^\infty (-1)^{i+1}i^2(4m-i^2)\overline{p}(4m-i^2)$$
when we take the sum $\pmod 8$. In this case, the sum is over terms of the form
$$(-1)^{i+1} 3 \overline{p}(4k+3)$$
and \eqref{overp1} follows by induction. 
\end{proof}
\begin{Remark}
Equation \eqref{overp1} follows from a theorem of Hirschhorn and Sellers
\cite[Equation (6)]{HS2005}. Thus \eqref{overp2} holds too. 
\end{Remark}
Theorems~\ref{th:equiv1} and \ref{th:equiv2} illustrate
the
 correspondence between results in the theory of partitions theory and results for divisor functions.


\begin{thebibliography}{10}

\bibitem{AAR1999}
G.~E. Andrews, R.~Askey, and R.~Roy.
\newblock {\em Special functions}, volume~71 of {\em Encyclopedia of
  Mathematics and its Applications}.
\newblock Cambridge University Press, Cambridge, 1999.

\bibitem{AE2004}
G.~E. Andrews and K.~Eriksson.
\newblock {\em Integer partitions}.
\newblock Cambridge University Press, Cambridge, 2004.

\bibitem{BB2022a}
H.~S. Bal and G.~Bhatnagar.
\newblock The {P}artition-{F}requency enumeration matrix.
\newblock {\em Ramanujan J. (published online)}, page 29pp, 2022.
\newblock
  \href{https://arxiv.org/abs/2102.04191}{https://arxiv.org/abs/2102.04191}.

\bibitem{Berndt1989}
B.~C. Berndt.
\newblock {\em Ramanujan's notebooks. {P}art {II}}.
\newblock Springer-Verlag, New York, 1989.

\bibitem{Berndt1994}
B.~C. Berndt.
\newblock {\em Ramanujan's notebooks. {P}art {IV}}.
\newblock Springer-Verlag, New York, 1994.

\bibitem{Berndt2006}
B.~C. Berndt.
\newblock {\em Number theory in the spirit of {R}amanujan}, volume~34 of {\em
  Student Mathematical Library}.
\newblock American Mathematical Society, Providence, RI, 2006.

\bibitem{bonciocat2003}
N.~C. Bonciocat.
\newblock Congruences for the convolution of divisor sum function.
\newblock {\em Bull. Greek Math. Soc.}, 47:19--29, 2003.

\bibitem{Christopher2019}
A.~D. Christopher.
\newblock Euler-type recurrence relation for arbitrary arithmetical function.
\newblock {\em Integers}, 19:Paper No. A62, 15, 2019.

\bibitem{CL2004}
S.~Corteel and J.~Lovejoy.
\newblock Overpartitions.
\newblock {\em Trans. Amer. Math. Soc.}, 356(4):1623--1635, 2004.

\bibitem{LE1760-243}
L.~Euler.
\newblock Observatio de summis divisorum.
\newblock {\em Novi Commentarii Academiae scientiarum Imperialis Petropoli-
  tanae}, 5:59--74, 1760.
\newblock Reprinted Euler Archive -- All Works E243, English translation (by
  Jordan Bell) available at \url{https://arxiv.org/abs/math/0411587v3}.

\bibitem{Fine1988}
N.~J. Fine.
\newblock {\em Basic hypergeometric series and applications}, volume~27 of {\em
  Mathematical Surveys and Monographs}.
\newblock American Mathematical Society, Providence, RI, 1988.
\newblock With a foreword by George E. Andrews.

\bibitem{gallardo2007}
L.~H. Gallardo.
\newblock On {B}onciocat's congruences involving the sum of divisors function.
\newblock {\em Bull. Greek Math. Soc.}, 53:69--70, 2007.

\bibitem{Glaisher1885}
J.~W.~L. {Glaisher}.
\newblock {On certain sums of products of quantities depending upon the
  divisors of a number.}
\newblock {\em {Mess.}}, (XV):1--20, 1885.

\bibitem{Glaisher1891}
J.~W.~L. {Glaisher}.
\newblock {Expressions for the sum of the cubes of the divisors of a number in
  terms of partitions of inferior numbers.}
\newblock {\em Mess.}, 2(XXI):47--48, 1891.

\bibitem{Gould1974}
H.~W. Gould.
\newblock Coefficient identities for powers of {T}aylor and {D}irichlet series.
\newblock {\em Amer. Math. Monthly}, 81:3--14, 1974.

\bibitem{Hahn2007}
H.~Hahn.
\newblock Convolution sums of some functions on divisors.
\newblock {\em Rocky Mountain J. Math.}, 37(5):1593--1622, 2007.

\bibitem{HS2005}
M.~D. Hirschhorn and J.~A. Sellers.
\newblock Arithmetic relations for overpartitions.
\newblock {\em J. Combin. Math. Combin. Comput.}, 53:65--73, 2005.

\bibitem{HOSW2002}
J.~G. Huard, Z.~M. Ou, B.~K. Spearman, and K.~S. Williams.
\newblock Elementary evaluation of certain convolution sums involving divisor
  functions.
\newblock In {\em Number theory for the millennium, {II} ({U}rbana, {IL},
  2000)}, pages 229--274. A K Peters, Natick, MA, 2002.

\bibitem{Melfi1998}
G.~Melfi.
\newblock On some modular identities.
\newblock In {\em Number theory ({E}ger, 1996)}, pages 371--382. de Gruyter,
  Berlin, 1998.

\bibitem{Merca2018a}
M.~Merca.
\newblock New connections between functions from additive and multiplicative
  number theory.
\newblock {\em Mediterr. J. Math.}, 15(2):Paper No. 36, 13, 2018.

\bibitem{Merca2021b}
M.~Merca.
\newblock Congruence identities involving sums of odd divisors function.
\newblock {\em Proc. Rom. Acad. Ser. A Math. Phys. Tech. Sci. Inf. Sci.},
  22(2):119--125, 2021.

\bibitem{Merca2021a}
M.~Merca.
\newblock Overpartitions and functions from multiplicative number theory.
\newblock {\em Politehn. Univ. Bucharest Sci. Bull. Ser. A Appl. Math. Phys.},
  83(3):97--106, 2021.

\bibitem{ORW1995}
K.~Ono, S.~Robins, and P.~T. Wahl.
\newblock On the representation of integers as sums of triangular numbers.
\newblock {\em Aequationes Math.}, 50(1-2):73--94, 1995.

\bibitem{Venkov1970}
B.~A. Venkov.
\newblock {\em Elementary number theory}.
\newblock Translated from the Russian and edited by Helen Alderson.
  Wolters-Noordhoff Publishing, Groningen, 1970.

\bibitem{Wilf2006}
H.~S. Wilf.
\newblock {\em generatingfunctionology}.
\newblock A K Peters, Ltd., Wellesley, MA, third edition, 2006.

\bibitem{Williams2011}
K.~S. Williams.
\newblock {\em Number theory in the spirit of {L}iouville}, volume~76 of {\em
  London Mathematical Society Student Texts}.
\newblock Cambridge University Press, Cambridge, 2011.

\end{thebibliography}
\end{document}